\documentclass[12pt]{amsart}
\usepackage{amssymb,amsthm,amsmath}
\usepackage{hyperref}
\usepackage{ifthen}
\usepackage{graphicx}
\nonstopmode
 \numberwithin{equation}{section}
\newtheorem{thm}{Theorem}[section]

\newtheorem{lem}[thm]{Lemma}

\newenvironment{pf}[1][]{%
	\vskip 3mm
	\noindent
	\ifthenelse{\equal{#1}{}}%
	{{\slshape Proof. }}%
	{{\slshape #1.} }%
}%
{\qed\bigskip}

\newtheorem{prop}[thm]{Proposition}

\allowdisplaybreaks[2]

\title[]
{Strichartz Inequality for orthonormal functions associated with special Hermite operator}

\author{Shyam Swarup Mondal}
\author{Jitendriya Swain}

\address{Shyam Swarup Mondal  \endgraf Department of Mathematics
	\endgraf IIT Guwahati
	\endgraf Guwahati, Assam, India.}
\email{mondalshyam055@gmail.com}

\address{Jitendriya Swain,  Associate professor  \endgraf Department of Mathematics
	\endgraf IIT Guwahati
	\endgraf Guwahati, Assam, India.}
\email{jitumath@iitg.ac.in}

\date{\today}
\begin{document}
	
	\maketitle
	\allowdisplaybreaks

	\begin{abstract} In this article, we obtain the Strichartz estimate  for the system of orthonormal functions associated with the special Hermite operator.
\end{abstract}
%	\tableofcontents 	
\section{Introduction}

Consider the free Schr\"odinger equation \begin{align}\label{502}
i \partial_{t} u(x, t)&=-\Delta u(x, t) \quad x \in \mathbb{R}^{n}, t \in \mathbb{R}\\\nonumber u(x, 0)&=f(x).
\end{align}
It is well known that $e^{i t \Delta} f$ is the unique solution to the initial value problem (\ref{502}).  The following remarkable estimate for the solution to the  initial value problem (\ref{502}) is first obtained
by Strichartz \cite{st} in connection with Fourier restriction theory:
\begin{thm}\label{V} Let $f\in L^2(\mathbb{R}^n)$. If  $p, q \geq 1$ satisfying $(p, q, n) \neq(1, \infty, 2)$ and
	$
	\frac{2}{p}+\frac{n}{q}=n,
	$ then $e^{i t \Delta} f\in L^{2p}_tL^{2q}_x(\mathbb{R}\times \mathbb{R}^n)$ and satisfies the inequality
	$$\int_{\mathbb{R}}\left(\int_{\mathbb{R}^{n}}\left|\left(e^{i t \Delta} f\right)(x)\right|^{2 q} d x\right)^{\frac{p}{q}} d t \leq  C\left(\int_{\mathbb{R}^{n}}|f(x)|^{2} d x\right)^{p}.
	$$\end{thm}

The above inequality have been substantially generalized  for  a system of orthonormal functions by Frank-Lewin-Lieb-Seiringer \cite{frank} and Frank-Sabin \cite{frank1}.
\begin{thm}\label{507}\cite{frank, frank1}
	Assume that
	$p, q, n \geq 1$ such that
	$$
	1\leq q < \frac{n+1}{n-1} \quad \text { and } \quad \frac{2}{p}+\frac{n}{q}=n.
	$$
	For any (possibly infinite) system $u_j$ of orthonormal functions in $L^{2}\left(\mathbb{R}^{n}\right)$
	and any coefficients $\left(n_{j}\right) \subset \mathbb{C},$ we have
	$$\int_{\mathbb{R}} \left( \int_{\mathbb{R}^{n}} \left|\sum_{j} n_{j}\left| \left(e^{i t \Delta} u_{j}\right) (x)  \right|^{2} \right|^qdx \right)^{\frac{p}{q}} d t \leq  C_{n, q}^{p}\left(\sum_{j}\left|n_{j}\right|^{\frac{2 q}{  q+1}}\right)^{\frac{p(q+1)}{2 q}}, $$
	where $C_{n, q}$ is a universal constant which only depends on $n$ and $q$.
\end{thm}
Further, Theorem \ref{V} has been extended for the Schr\"odinger equation associated with the Hermite operator $H=-\Delta+|x|^2$:
\begin{align}\label{503}
i \partial_{t} u(x, t)&=H u(x, t) \quad x \in \mathbb{R}^{n}, t \in \mathbb{R}\\\nonumber u(x, 0)&=f(x).
\end{align}
Assuming $f\in L^2(\mathbb{R}^n)$, the solution of the initial value problem (\ref{503}) is given  by $u(x,t)=e^{-i t H} f(x).$ The Strichartz inequality in this case has been  proved by Koch-Tataru \cite{KT} or Nandakumaran-Ratnakumar \cite{Ratna} resulting in the following.
\begin{thm} \label{504}
	Let $f\in L^2(\mathbb{R}^n)$. If  $p, q \geq 1$ satisfying $(p, q, n) \neq(1, \infty, 2)$ and
	$\frac{2}{p}+\frac{n}{q}=n,$ then $$\|e^{-i t H} f\|_{L_{t}^{2p} L_{x}^{2q}\left(\mathbb{T}\times \mathbb{R}^{n}\right)}\leq C\|f\|_{2}.$$
\end{thm}
The above inequality  is further generalized for  a system of orthonormal functions in \cite{lee,ssm}.
\begin{thm}\cite{lee,ssm}\label{STH}[Strichartz inequality for orthonormal functions for Hermite operator]
	Let
	$p, q, n \geq1$ such that  $$ 1\leq q < \frac{n+1}{n-1}\quad \text { and } \quad \frac{2}{p}+\frac{n}{q}=n.$$ For any (possibly infinite) system $\left(u_{j}\right)$ of orthonormal functions in $L^{2}\left(\mathbb{R}^{n}\right)$
	and any coefficients $ \left( n_{j} \right) \subset \mathbb{C},$ we have
	\begin{align}\label{2}
	\int_{-\pi}^{\pi}\left( \int_{\mathbb{R}^{n}}\left| \sum_{j} n_{j}\left| \left(e^{-i t H} u_{j}\right) ( x )\right|^{2} \right|^{q} d x \right)^{\frac{p}{q}} d t  \leq C_{n, q}^{p}\left(\sum_{j}\left|n_{j}\right|^{\frac{2 q}{q+1}}\right)^{\frac{p(q+1)}{2 q}},
	\end{align}
	where $C_{n, q}$ is a universal constant which only depends on $n$ and $q$.
\end{thm}

 The Strichartz estimate for the Schr\"odinger equation associated with special Hermite operator  $ \mathcal{L}$ (defined in section 2)
 on $L^2(\mathbb{C}^{n})$  has been considered by Ratnakumar \cite{Ratna111} in the following initial value problem:
\begin{align}\label{s1}
i \partial_{t} u(z, t)&=\mathcal{L}u(z, t) \quad z \in \mathbb{C}^{n}, t \in \mathbb{R}\\\nonumber u(z, 0)&=f(z).
\end{align}
For  $f\in L^2(\mathbb{C}^n)$, the solution of the initial value problem (\ref{s1}) is given  by $u(z,t)=e^{-i t \mathcal{L}} f(z)$ and satisfies the following Strichartz estimate.
\begin{thm}\cite{Ratna111}\label{504}
	Let $f\in L^2(\mathbb{C}^n)$.  	If $1<p<\infty,~ \frac{1}{p} \geq n\left(1-\frac{1}{q}\right)$$~\text{or }~
\frac{1}{2} \leq p \leq 1, ~1\leq q<\frac{n}{n-1}$ then
	$$\|e^{-i t \mathcal{L}} f\|_{L_{t}^{2p} L_{z}^{2q}\left(\mathbb{T}\times \mathbb{C}^{n}\right)}\leq C\|f\|_{2}.$$
\end{thm}

The main aim of this paper is to obtain the following Strichartz estimate for system of orthonormal functions associated with the special operator with respect to the special Hermite transform. To the best of our knowledge the study on restriction theorem with respect to the special Hermite transform has not been considered in the literature so far.

\begin{thm}\label{STHH}
	Let $ q, n \geq 1$ and $p>1$ such that  $$ 1\leq q \leq  1+\frac{1}{n}\quad \text { and } \quad \frac{1}{p}+\frac{n}{q}=n.$$ For any (possibly infinite) system $\left(u_{j}\right)$ of orthonormal functions in $L^{2}\left(\mathbb{C}^{n}\right)$ and any coefficients $\left(n_{j}\right) \subset \mathbb{C},$ there exists a constant $C>0$ such that
\begin{align}\label{SH1}
\left\|\sum_{j} n_{j}\left|e^{-i t \mathcal{L} } u_{j}\right|^{2}\right\|_{L_{t}^{p} L_{z}^{q}\left(\mathbb{T}\times \mathbb{C}^{n}\right)} \leq C \left(\sum_{j}\left|n_{j}\right|^{\frac{2q}{q+1}}\right)^{\frac{(q+1)}{2 q}}.
	\end{align}
	\end{thm}
Let $f\in L^1(\mathbb{C}^n)$. Define the special Hermite  transform of $f$ by $$\hat{f}(\mu, \nu)=\int_{\mathbb{C}^n}f(z)\Phi_{\mu\nu}(z)\,dz, \quad \mu , \nu\in \mathbb{N}_0^n$$ where   $\mathbb{N}_0$ denotes the set of all non-negative integers and  $\Phi_{\mu\nu}$'s are the special Hermite functions (defined in section 2) on $\mathbb{C}^n$.
If $f\in L^2(\mathbb{C}^n)$  then $\{\hat{f}(\mu, \nu)\}\in \ell^2(\mathbb{N}_0^{2n})$ and satisfies the Plancherel formula  $$\|f\|_2^2=\sum_{(\mu,\nu)\in{\mathbb{N}^{2n}_0}}|\hat{f}(\mu, \nu)|^2.$$  The inverse special Hermite transform is given by $$f(z)=\sum_{(\mu, \nu) \in\mathbb{N}^{2n}_0}\hat{f}(\mu, \nu)\Phi_{\mu\nu}(z).$$  Given a discrete surface $S$ in $\mathbb{N}_0^{2n}$, we define the restriction operator $(\mathcal{R}_Sf):=\{\hat{f}(\mu, \nu)\}_{\mu, \nu\in S}$  and the operator dual to $\mathcal{R}_S$ (called the extension operator)  as
 $$\mathcal{E}_S(\{\hat{f}(\mu, \nu)\}):= \sum_{\mu, \nu\in S} \hat{f}(\mu, \nu)\Phi_{\mu\nu}.$$
  We consider the following problem:

 {\bf Problem 1:} For which exponents $1\leq p\leq 2,$ the sequence of special Hermite transforms of a function $f\in L^p(\mathbb{C}^n)$ belongs to $\ell^2(S)$?

 This question can be reframed to the boundedness of the operator $\mathcal{E}_S$ from $\ell^2(S)$ to $L^{p'}(\mathbb{C}^n)$, where $p'$ is the conjugate exponent of $p$. Since $\mathcal{E}_S$ is bounded from $\ell^2(S)$ to $L^{p'}(\mathbb{C}^n)$ if and only if $T_S:=\mathcal{E}_S(\mathcal{E}_S)^*$ is bounded from $L^p(\mathbb{C}^n)$ to $L^{p'}(\mathbb{C}^n)$,
 Problem 1 can be re-written as follows:

 {\bf Problem 2:} For which exponents $1\leq p\leq 2,$ the operator  $T_S:=\mathcal{E}_S(\mathcal{E}_S)^*$ is bounded from $L^p(\mathbb{C}^n)$ to $L^{p'}(\mathbb{C}^n)$?

 To address this problem we introduce an analytic family of operators $(T_z)$ defined on the strip $a\leq {\rm Re}~ z\leq b$ in the complex plane  such that $T_S=T_c$ for some $c\in (a, b)$ and show that the operator $W_1T_SW_2$ belongs to a Schatten class for $W_1, W_2\in L^{\frac{2p}{2-p}}(\mathbb{T}\times\mathbb{C}^n)$, which is more general  result $L^p-L^{p'}$ boundedness of  $T_S.$

 Such problems are often considered in the literature. For example, on $\mathbb{R}^n$, the celebrated Stein-Tomas Theorem (see \cite{Stein,T1,T2}) gives an affirmative answer to Fourier restriction problem for compact surfaces with non-zero Gaussian curvature if and only if $1\leq p\leq \frac{2(n+1)}{n+3}$.  For quadratic surfaces, Strichartz \cite{st}  gave a complete solution to Fourier restriction problem, when $S$ is a quadratic surface given by $S = \{x \in \mathbb{R}^n: R(x) = r\},$ where $R(x)$ is a polynomial of
 degree two with real coefficients and $r$ is a real constant. Further the Stein-Tomas Theorem is generalized to a system of orthonormal functions with respect to the Fourier transform by Frank-Lewin-Lieb-Seiringer \cite{frank1} and Frank-Sabin \cite{frank}.

The schema of the paper apart from introduction  is as follows: In Section 2,  we  discuss the spectral theory of the  Hermite operator and the kernel estimates for the special Hermite semigroup. In section 3, we obtain the duality principle in terms of Schatten bounds of the operator $We^{-it\mathcal{L}}(e^{-it\mathcal{L}})^*\overline{W}$
and prove the Strichartz estimate for $1\leq q\leq1+\frac{1}{n}$, for the system of  orthonormal functions associated with  the special Hermite operator as the restriction of the special Hermite  transform to the discrete surface $S=\{(\mu, \nu, \lambda)\in \mathbb{N}_0^n\times\mathbb{N}_0^n \times \mathbb{Z}: \lambda=2|\nu|+n\}$.
\section{Preliminary}
In this section we discuss some basic definitions and provide necessary background information  about the special Hermite semigroup.
\subsection{Hermite Operator and Special Hermite functions}
Let $\mathbb{N}_0$ be the set of all non-negative integers. Let $H_k$ denote the Hermite polynomial on $\mathbb{R}$, defined by
$$H_k(x)=(-1)^k \frac{d^k}{dx^k}(e^{-x^2} )e^{x^2}, \quad k\in \mathbb{N}_0$$
and $h_k$ denote the normalized Hermite functions on $\mathbb{R}$ defined by
$$h_k(x)=(2^k\sqrt{\pi} k!)^{-\frac{1}{2}} H_k(x)e^{-\frac{1}{2}x^2}, \quad k\in \mathbb{N}_0.$$ %The Hermite functions $\{h_k \}$ are the eigenfunctions of the Hermite operator $H=-\frac{d^2}{dx^2}+x^2$ with eigenvalues $2k+1,  k=0, 1, 2, \cdots$. These functions form an orthonormal basis for $L^2(\mathbb{R})$.
The higher dimensional Hermite functions denoted by $\Phi_{\alpha}$ are obtained by taking tensor product of one dimensional Hermite functions. Thus for any multi-index $\alpha \in \mathbb{N}_0^n$ and $x \in \mathbb{R}^n$, we define
 $\Phi_{\alpha}(x)=\prod_{j=1}^{n}h_{\alpha_j}(x_j).$

For each multi-index $\mu,\nu$ and $\zeta\in\mathbb{C},$ we define the special Hermite functions $\Phi_{\mu \nu}$ by
$$\Phi_{\mu \nu}(\zeta)=(2 \pi)^{-\frac{n}{2}} \int_{\mathbb{R}^{n}} e^{i x \cdot \xi} \Phi_\mu\left(\xi+\frac{y}{2}\right)  \Phi_\nu\left(\xi-\frac{y}{2}\right) d \xi, ~\zeta=x+i y \in \mathbb{C}^{n}.$$ The family of functions $\{\Phi_{\mu\nu}\}$ form an orthonormal basis for $L^2(\mathbb{C}^n).$ The special Hermite functions are the eigenfunctions of the special Hermite operator $\mathcal{L}$ (or the twisted Laplacian)  defined by
$$
\mathcal{L}=\frac{1}{2} \sum_{j=1}^{n}\left(Z_{j} \bar{Z}_{j}+\bar{Z}_{j} Z_{j}\right),
$$ where $Z_{j}=\frac{\partial}{\partial \zeta_{j}}+\frac{1}{2} \bar{\zeta}_{j}, \bar{Z}_{j}=-\frac{\partial}{\partial \bar{\zeta}_{j}}+\frac{1}{2} \zeta_{j}, j=1,2, \ldots n$  with eigenvalues $(2|\nu|+n)$.
The special Hermite operator $\mathcal{L}$ is
self-adjoint and admits a spectral decomposition in terms of special
Hermite functions. Given \(f \in L^{2}\left(\mathbb{C}^{n}\right)\) the expansion
\begin{align}\label{sh2}
	f=\sum_{\mu ,\nu\in \mathbb{N}^{n}_0}\langle f, \Phi_{\mu \nu}\rangle  \Phi_{\mu \nu}
\end{align}
%converges to \(f\) in \(L^{2}\left(\mathbb{C}^{n}\right) .\) For each \(k \in \mathbb{N}\), defining \(P_{k}\) to be the
%orthogonal projection of \(L^{2}\left(\mathbb{C}^{n}\right)\) onto the eigenspace spanned by \(\left\{\Phi_{\mu \nu} :|\nu|=k\right\}\),
%the spectral decomposition of \(\mathcal{L}\) can be written as
%$$
%\mathcal{L} f=\sum_{k=0}^{\infty}(2 k+n) P_{k} f.
%$$
%The twisted convolution of two functions $f$ and $g$ on $\mathbb{C}^{n}$ is defined by
%$$
%f \times g(\zeta)=\int_{\mathbb{C}^{n}} f(\zeta-w) g(w) e^{\frac{i}{2} \operatorname{Im}(\zeta \cdot \bar{w})} d w,
%$$
%where $\operatorname{Im}$ denotes the imaginary part. The family $\{\Phi_{\mu \nu} \}$ satisfies the following orthogonality properties
%\begin{align}\label{sh1}
%\Phi_{\mu \nu} \times \Phi_{\alpha \beta}=\left\{\begin{array}{ll}{(2 \pi)^{n / 2} \Phi_{\mu \beta} } & {\text{if } \nu=\alpha,} \\ {0} & {\text{otherwise}.} \end{array}\right.
%\end{align}
%Using these properties we have $$f(\zeta)=(2 \pi)^{-n} \sum_{k} f \times \phi_{k}(\zeta),$$ where $\phi_{k}(\zeta)=(2 \pi)^{n / 2} \displaystyle\sum_{|\nu|=k} \Phi_{\nu \nu}(\zeta)$.
converges to \(f\) in \(L^{2}\left(\mathbb{C}^{n}\right) .\)  The above expansion also can be   written as $f= \sum_{k=0}^\infty P_kf$, where
$$P_k=\sum_{\mu ,|\nu|=k}\langle \cdot , \Phi_{\mu \nu}\rangle  \Phi_{\mu \nu}$$
is the orthogonal projection of \(L^{2}\left(\mathbb{C}^{n}\right)\) onto the eigenspace spanned by \(\left\{\Phi_{\mu \nu} :|\nu|=k\right\}\).
%For each \(k \in \mathbb{N}\), defining \(P_{k}\) to be the
%orthogonal projection of \(L^{2}\left(\mathbb{C}^{n}\right)\) onto the eigenspace spanned by \(\left\{\Phi_{\mu \nu} :|\nu|=k\right\}\),
%the spectral decomposition of \(\mathcal{L}\) can be written as
For each \(k \in \mathbb{N}\),  the spectral decomposition of \(\mathcal{L}\) can be written as
$$
\mathcal{L} f=\sum_{k=0}^{\infty}(2 k+n) P_{k} f.
$$
The twisted convolution of two functions $f$ and $g$ on $\mathbb{C}^{n}$ is defined by
$$
f \times g(\zeta)=\int_{\mathbb{C}^{n}} f(\zeta-w) g(w) e^{\frac{i}{2} \operatorname{Im}(\zeta \cdot \bar{w})} d w,
$$
where $\operatorname{Im}$ denotes the imaginary part. The family $\{\Phi_{\mu \nu} \}$ satisfies the following orthogonality properties
\begin{align}\label{CH3sh1}
	\Phi_{\mu \nu} \times \Phi_{\alpha \beta}=\left\{\begin{array}{ll}{(2 \pi)^{n / 2} \Phi_{\mu \beta} } & {\text{if } \nu=\alpha,} \\ {0} & {\text{otherwise}.} \end{array}\right.
\end{align}
Let  $L_{k}^{\alpha}$ denote the Laguerre polynomial of degree $k$ and of  order $\alpha>-1,$ defined by  the generating function identity (see \cite{Ratna111})
$$
\sum_{k=0}^{\infty} L_{k}^{\alpha}(x) \omega^{k}=(1-\omega)^{-\alpha-1} e^{-\frac{\omega}{1-\omega} x}, \quad|\omega|<1
$$
and let   $\phi_{k}(z)=L_{k}^{n-1}\left(\frac{1}{2}|z|^{2}\right) e^{-\frac{1}{4}|z|^{2}}$ be the Laguerre function of order $n-1 .$
The special Hermite functions $\Phi_{\nu \nu}$  are related to the Laguerre functions $\phi_{k}$  by the following relation
\begin{align}\label{CH33333}
	(2 \pi)^{n / 2} \sum_{|\nu|=k} \Phi_{\nu \nu}=\phi_{k}.
\end{align}
Now taking twisted convolution on both sides of (\ref{sh2}) with $\Phi_{\alpha \alpha}$  and using the orthogonality property (\ref{CH3sh1}),  we have
\begin{align}\label{CH3333}
	f \times \Phi_{\alpha \alpha}=(2 \pi)^{n / 2} \sum_{\mu}\left\langle f, \Phi_{\mu \alpha}\right\rangle \Phi_{\mu \alpha}.
\end{align}
Summing both sides of (\ref{CH3333}) with respect to all $\alpha$ such that $|\alpha|=k$ and using (\ref{CH33333}),   the spectral projection $P_{k}$  has the simpler representation
$$
P_{k} f(\zeta)=(2 \pi)^{-\frac{n}{2}} \sum_{|\alpha|=k} f \times \Phi_{\alpha \alpha}(\zeta)=(2 \pi)^{-n} f \times \varphi_{k}(\zeta).
$$
Then the special Hermite expansion takes the compact form $$f(\zeta)=(2 \pi)^{-n} \sum_{k} f \times \phi_{k}(\zeta).$$
The operator $\mathcal{L}$ defines a semigroup, called the special Hermite semigroup and denoted by
$e^{-t \mathcal{L}}, t>0$,   by the expansion
$$
e^{-t \mathcal{L}} f=(2 \pi)^{-\frac{n}{2}}\sum_{k=0}^{\infty} e^{-(2 k+n) t} f \times \phi_{k}
$$
for $f \in L^{2}(\mathbb{C}^{n}).$
For the auxiliary complex semigroup  $\{e^{-\eta \mathcal{L}}\}$, $\eta=r+i t, ~r>0$, we write
$$
e^{-\eta \mathcal{L}} f(\zeta)=(2 \pi)^{-n} \sum_{k=0}^{\infty} e^{-\eta(2 k+n)} f \times \phi_{k}(\zeta).
$$
Thus,  $e^{-\eta \mathcal{L}}$  is a  twisted convolution operator
$$
e^{-\eta \mathcal{L}} f(\zeta)=\int_{\mathbb{C}^{n}} f(\zeta-w) K_{\eta}(w) e^{\frac{i}{2} \operatorname{Im} (\zeta \cdot \bar{w})} d w
$$
with kernel (see \cite{Ratna111})
$$
K_{\eta}(\zeta)=(2 \pi)^{-n} \sum_{k=0}^{\infty} e^{-\eta(2 k+n)} \phi_{k}(\zeta)=(2 \pi)^{-n} e^{-n \eta}(1-\omega)^{-n} e^{-\frac{1+\omega}{1-\omega} \frac{|\zeta|^{2}}{4}},
$$
where $\omega=e^{-2 \eta}.$ So  $K_{r+it}(\zeta)=K_{r+i(t+2\pi)}(\zeta)$, and
\begin{align}\label{kk}
\left|K_{\eta}(\zeta)\right| \leq \frac{2}{|\sin t|^{n}}, \quad \eta=r+i t, \quad \zeta \in \mathbb{C}^{n}.
\end{align}
We refer to \cite{Ratna111} for a detailed study on special Hermite semigroup.

\subsection{Schatten class and the duality principle} Let $\mathcal{H}$ be a complex and separable Hilbert space in which the inner product is denoted by $\langle, \rangle_\mathcal{H}$. Let $T:\mathcal{H} \rightarrow \mathcal{H}$ be a compact operator and let    $T^{*}$ denotes the adjoint of $T$.  For $1 \leq r<\infty,$ the Schatten space $\mathcal{G}^{r}(\mathcal{H})$ is defined as the space of all compact operators $T$ on $\mathcal{H}$ such that $$\sum_{n =1}^\infty  \left(s_{n}(T)\right)^{r}<\infty,$$ where $s_{n}(T)$ denotes the singular values of \(T,\) i.e., the eigenvalues of \(|T|=\sqrt{T^{*} T}\) counted according to multiplicity. For $T\in \mathcal{G}^{r}(\mathcal{H})$, the Schatten $r$-norm is defined by $$ \|T\|_{\mathcal{G}^{r}}=\left(\sum_{n=1}^{\infty}\left(s_{n}(T)\right)^{r}\right)^{\frac{1}{r}}. $$An operator belongs to the class \(\mathcal{G}^{1}({\mathcal{H}})\) is known as {\it Trace class} operator. Also, an operator belongs to   \(\mathcal{G}^{2}({\mathcal{H}})\) is known as  {\it Hilbert-Schmidt} operator.

 \section{Strichartz inequality for system of  orthonormal functions}
 In order to obtain the Strichartz inequality for the system of orthonormal functions we need a duality principle lemma in our context. We refer to Proposition 1 and Lemma 3 of \cite{frank1} with appropriate modifications to obtain the following two results:

 \begin{prop}\label{sh4}
 	Let $(T_z)$ be an analytic family of operators on $ \mathbb{T}\times \mathbb{C}^{n}$ in the sense of Stein defined on the strip $-\lambda_0\leq  \operatorname{Re} z\leq 0$
 	for some $\lambda_0 > 1$. Assume that we have the following bounds
 	\begin{align}\label{expo}
 	&\left\|T_{i s}\right\|_{L^{2}(\mathbb{T}\times \mathbb{C}^n ) \rightarrow L^{2}(\mathbb{T}\times \mathbb{C}^n )} \leq M_{0} e^{a|s|},\\\nonumber\\\label{expo1}
 	&\left\|T_{-\lambda_{0}+i s}\right\|_{L^{1}(\mathbb{T}\times \mathbb{C}^n ) \rightarrow L^{\infty}(\mathbb{T}\times \mathbb{C}^n )} \leq M_{1} e^{b|s|}
 	\end{align}
 	for all $s \in \mathbb{R}$ and 	for some $a, b, M_{0}, M_{1} \geq 0 $. Then, for all  $W_{1}, W_{2} \in L^{2 \lambda_{0}}\left(\mathbb{T}\times \mathbb{C}^n , \mathbb{C}\right)$
 	the operator $W_{1} T_{-1} W_{2}$ belongs to $\mathcal{G}^{2 \lambda_{0}}\left(L^{2}\left(\mathbb{T}\times \mathbb{C}^n  \right)\right)$ and we have the estimate
 	\begin{align}\label{510}
 	\left\|W_{1} T_{-1} W_{2}\right\|_{\mathcal{G}^{2 \lambda_{0}}\left(L^{2}\left(\mathbb{T}\times \mathbb{C}^n \right)\right)} \leq  M_{0}^{1-\frac{1}{\lambda_{0}}} M_{1}^{\frac{1}{\lambda_{0}}}\left\|W_{1}\right\|_{L_t^{2 \lambda_{0}}L_z^{2 \lambda_{0}}\left(\mathbb{T}\times \mathbb{C}^n \right)}\left\|W_{2 }\right\|_{L_t^{2 \lambda_{0}}L_z^{2 \lambda_{0}}\left(\mathbb{T}\times \mathbb{C}^n \right)}.
 	\end{align}
 \end{prop}

%We state the following duality principle which can be obtained with a trivial modification to the duality principle due to Frank-Sabin \cite{frank1}.

\begin{lem}\label{dual}(Duality principle)
	Let $p,q\geq 1$ and $\alpha\geq1$. Let $Af(t,w)=e^{-it\mathcal{L}}f(w)$. Then the following statements are equivalent.
	\begin{enumerate}
		\item There is a constant \(C>0\) such that
		\begin{align}\label{511}
		\left\|W A A^{*} \overline{W}\right\|_{\mathcal{G}^{\alpha}\left(L^{2}\left( \mathbb{T} \times\mathbb{C}^{n} \right)\right)} \leq C\|W\|_{{L_t^{\frac{2q}{2-q}}L_w^{\frac{2p}{2-p}}(\mathbb{T}\times \mathbb{C}^{n})}}^{2} \end{align} for all $W \in {L_t^{\frac{2q}{2-q}}L_w^{\frac{2p}{2-p}}(\mathbb{T}\times \mathbb{C}^{n})}$, where the function $W$ is interpreted as an operator which acts by multiplication.
		
		\item  For any orthonormal system $\left(f_{j}\right)_{j \in J}$
		in  $L^2(\mathbb{C}^n)$ and any sequence $\left(n_{j}\right)_{j \in J} \subset \mathbb{C}$, there is a constant \(C'>0\) such that
		\begin{align}\label{512}
		\left\|\sum_{j \in J} n_{j} \left| A f_{j}\right|^{2}\right\|_{L_t^{\frac{q'}{2}}L_w^{\frac{p'}{2}}(\mathbb{T}\times \mathbb{C}^{n})} \leq C' \left(\sum_{j \in J}\left|n_{j}\right|^{\alpha^{\prime}}\right)^{1 / \alpha^{\prime}}.
		\end{align}
	\end{enumerate}
	\end{lem}
 If  $f\in L^2(\mathbb{C}^n)$, the solution of the initial value problem (\ref{s1}) can be realized as the extension operator of some function $f$ on $\mathbb{T}\times\mathbb{C}^n.$

Let $S$ be the discrete surface  $S=\{(\mu, \nu, \lambda)\in \mathbb{N}_0^n\times \mathbb{N}_0^n\times \mathbb{Z}: \lambda=2|\nu|+n\}$ with respect to counting measure. Then for all $f$ such that $\hat{f} \in \ell^{1}(S)$ and for all $(t, \zeta)\in [-\pi,\pi]\times \mathbb{C}^n$, the extension operator can be written as
\begin{align}\label{surfacee}
\mathcal{E}_{S} f(t, \zeta)=\sum_{\mu, \nu, \lambda\in S}\hat{f}(\mu, \nu, \lambda) \Phi_{\mu\nu} (\zeta)e^{-it\lambda},
\end{align} where $\hat{f}(\mu, \nu, \lambda)=\int_{\mathbb{C}^n}\int_{\mathbb{T}}f(t, w)\Phi_{\mu\nu}(w)e^{i\lambda t}\,dt dw $.
Using the fact that $$f\times\Phi_{\mu\mu}=(2\pi)^{\frac{n}{2}}\sum_{\nu}\langle f, \Phi_{\mu \nu}\rangle  \Phi_{\mu \nu}(\zeta)$$ and choosing
$$
\hat{f}(\mu, \nu, \lambda)=\left\{\begin{array}{ll}       {(2\pi)^n\hat{u}(\mu, \nu)} & {\text { if } \lambda=2|\nu|+n,}\\{0} & {~~\text {otherwise} }\end{array} \right.
$$ for some $u:\mathbb{C}^n\to \mathbb{C}$
in (\ref{surfacee}), we get
\begin{align*}
\mathcal{E}_{S} f(t, \zeta)&=(2\pi)^n\sum_{\mu, \nu\in S}\hat{u}(\mu, \nu) \Phi_{\mu\nu}(\zeta)e^{-it (2|\nu|+n)}\\
&=(2\pi)^n\sum_{\nu} \left( \sum_{\mu}\langle u, \Phi_{\mu \nu}\rangle  \Phi_{\mu \nu}(\zeta)\right) e^{-it (2|\nu|+n)}\\
&=(2\pi)^{\frac{n}{2}}\sum_{\nu}  e^{-it (2|\nu|+n)}u\times  \Phi_{\nu \nu}(\zeta)~\\
%&=(2\pi)^{-\frac{n}{2}}\sum_{k=0}^{\infty} e^{-it (2k+n)}\sum_{|\nu|=k} f\times  \phi_{\nu \nu}(z)~\\
&=(2\pi)^{\frac{n}{2}}\sum_{k=0}^{\infty} e^{-it (2k+n)}\left( u\times  \sum_{|\nu|=k} \Phi_{\nu \nu}(\zeta)\right)\\
&=\sum_{k=0}^{\infty} e^{-it (2k+n)} u\times   \phi_k(z)=e^{-it\mathcal{L}}u(\zeta).
\end{align*}
For $-1< {\rm Re}~ z\leq 0$, define the analytic family of generalized functions
\begin{align}\label{s11}
G_z(\mu , \nu, \lambda )=\frac{1}{\Gamma(z+1)}(\lambda-(2|\nu|+n))_+^z,
\end{align} where $$(\lambda-(2|\nu|+n))_+^z=\begin{cases}
 	(\lambda-(2|\nu|+n))^z~\text{ ~for }  \lambda-(2|\nu|+n)>0,\\0~\quad\quad\quad \quad\quad\qquad\text{ for }\lambda-(2|\nu|+n)\leq0.
 	\end{cases} $$

For Schwartz class functions $\phi$ on $\mathbb{N}_0^{2n}\times\mathbb{Z},$ using the discrete Taylor series expansion (see \cite{T}), we have
 \begin{align}
 	&	\left\langle \varsigma_{+}^{z}, \phi\right\rangle\label{CH31004}=\sum_{\varsigma\in \mathbb{N}_0^{2n}\times\mathbb{Z}}  \varsigma_{+}^{z}\phi(\varsigma)\\
 	&\label{CH31001}=\sum_{\varsigma\in \mathbb{N}_0^{2n}\times\mathbb{Z}}  \varsigma_{+}^{z}\left[ \phi(\varsigma)- \sum_{|\alpha|<M}\frac{1}{\alpha!}\varsigma^\alpha\Delta^\alpha\phi(0)                \right] +\sum_{|\alpha|<M}\frac{1}{\alpha!}\Delta^\alpha\phi(0)\sum_{\varsigma\in \mathbb{N}_0^{2n}\times\mathbb{Z}}  \varsigma^\alpha\varsigma_{+}^{z}.
 \end{align}
 The above formula is valid for $ z\neq-1, -2, \cdots,$
 regularizing (\ref{CH31004}).  Notice that (\ref{CH31001}) shows that $\left\langle \varsigma_{+}^{z}, \phi\right\rangle$ is treated as a function of $z$ with simple poles at $z=-1, -2, \cdots.$

 Thus for Schwartz class functions $\phi$ on $\mathbb{N}_0^{2n}\times\mathbb{Z},$ we have $$\displaystyle\lim_{z\to -1}\left\langle G_{z}, \phi\right\rangle=\lim_{z\to -1}\frac{1}{\Gamma(z+1)} \sum_{\mu,\nu}\phi(\mu, \nu, \lambda)(\lambda-(2|\nu|+n))_+^z=\sum_{(\mu,\nu, \lambda )\in S}\phi(\mu, \nu, \lambda).$$ We refer to \cite{sh} for the distributional calculus of $(\lambda-(2|\nu|+n))_+^z$.
Thus $G_{-1}=\delta_S$.
%\begin{align}\label{T_z}T_z g(x,y,  t) =  \sum_{\mu, \nu, \lambda} \hat{g}(\mu, \nu, \lambda) G_z( \mu, \nu, \lambda) \phi_{\mu\nu}(x, y) e^{-i\lambda t }.\end{align}Then we have
To prove our main result we need to prove the following proposition.
\begin{prop}\label{515}
Let $-1< {\rm Re}~ z<0$. Then the series
$ \sum_{k=0 }^\infty k_+^z e^{-i tk }$ is the Fourier series of a integrable function on $[-\pi, \pi]$ which is of class $C^\infty $ on $[-\pi, \pi]\setminus  \{0\}.$ Near  origin this function has the same singularity as the function whose values are $ \Gamma(z+1)(it)^{-z-1}$, i.e.,
\begin{align}\label{eq4} \sum_{k=0 }^\infty k_+^z e^{-i tk }\sim  \; \Gamma(z+1)(it)^{-z-1}+b(t),\end{align}
where $b\in C^\infty[-\pi , \pi]$.
 %Then\begin{align}\label{eq4}
%  \sum_{k=0 }^\infty k_+^z e^{-i tk }\asymp \; \Gamma(z+1)(it)^{-z-1}.
%\end{align}
 \end{prop}
\begin{pf}
 %let  us consider the expression
 %$\sum_{k=0 }^\infty k_+^z e^{-\tau k }e^{-itk}=\sum_{k=0 }^\infty k_+^z e^{-i s k },$ where $s=t-i\tau, $ and  $\tau= \operatorname{Im} (-s)$ so that $-\pi<\arg s<0$. Then the sum $\sum_{k=0 }^\infty k_+^z e^{-i s k }$ converges to  $\sum_{k=0 }^\infty k_+^z e^{-i t k }$  in the sense of distributions as$\tau \to 0$.
For $\tau>0$,  we calculate the inverse Fourier  transform of  $u_+^z e^{-\tau u }.$
  \begin{align*}
  \mathcal{F}^{-1}[ u_+^z e^{-\tau u }](x)&=\int_\mathbb{R}  u_+^z e^{-\tau u } e^{-i ux }\;du=\int_0^\infty   u^z e^{-i s u}\;du,
  \end{align*}
    where $s=x-i\tau$ so that $-\pi<\arg s<0$. Then   $ u_+^z e^{-\tau u }$ converges to  $u_+^z  $  in the sense of distributions as $\tau \to 0$. Also,  the inverse Fourier  transform of  $u_+^z e^{-\tau u }$ converges to  the inverse Fourier  transform of  $u_+^z  $. Using the change of variable $isu=\xi$ and proceeding as in page 170 of \cite{sh}, we get

    \begin{align}\label{tau}
  \mathcal{F}^{-1}[ u_+^z e^{-\tau u }](x)& =\frac{1}{(is)^{z+1}}\int_L  \xi^z e^{-\xi} \;d\xi=\frac{\Gamma(z+1)}{(is)^{z+1}},
\end{align}
  where   the contour $L$ of the integral is a ray from origin to infinity whose angle with respect to the real axis is given by $\arg\xi=\arg s+\frac{\pi}{2}$.
 Letting  $\tau \to 0$ in (\ref{tau}),  we have
   \begin{align}\label{eq50}
  \mathcal{F}^{-1}[ u_+^z ](x)& =\frac{\Gamma(z+1)}{(ix)^{z+1}}.
  \end{align}
We use the idea given in Theorem 2.17 of \cite{stein} to prove (\ref{eq4}). To make the paper self contained, we will only indicate the main steps. Let us consider a function $\eta\in C^\infty(\mathbb{R})$ such that $\eta(x)=1$ if $|x|\geq 1$, and vanishes in a neighborhood of the origin. Let $F(x)=\eta(x)x_{+}^z $ for $x\in  \mathbb{R}$. Writing  $F(x)=x_{+}^z +(\eta(x)-1)x_{+}^z $, using (\ref{eq50}) and denoting $f$ to be the inverse Fourier transform of $F$ in the sense of distributions, we have  $f(x)= {\Gamma(z+1)}{(ix)^{-z-1}}+b_1(x)$, where $b_1$ is the inverse Fourier transform of the   integrable function %whose values are
$(\eta(x)-1)x_{+}^z$  whose support is bounded.  Moreover, $b_1\in C^\infty(\mathbb{R})$ and $f\in L^1(\mathbb{R}).$

Applying Poisson summation formula (see page  250 of  \cite{stein}) to the function $f$ and using the fact   $\hat{f}=F$, we get
\begin{align*}
  \sum_{k=0 }^\infty k_+^z e^{-i tk }&= \sum_{k \in \mathbb{Z}} F(k)e^{-i tk }\\
  &\sim  \sum_{k\in \mathbb{Z}} f(2k\pi+t)  \\
  &=f(t)+ \sum_{|k|>0} f(2k\pi+t)  \\
  &= {\Gamma(z+1)}{(it)^{-z-1}}+b_1(t)+ \sum_{|k|>0} f(2k\pi+t)\\
  &= \Gamma(z+1)(it)^{-z-1}+b(t),
\end{align*}
where $b(t)=b_1(t)+ \sum_{|k|>0} f(2k\pi+t)\in C^\infty[-\pi , \pi]$.

\end{pf}

Now we are in a position to prove the following Strichartz inequality for the diagonal case.
\begin{thm}\label{diagonalcase}
	Let $n \geq 1$.  For any (possibly infinite) system $\left(u_{j}\right)$ of orthonormal functions in $L^{2}\left(\mathbb{C}^{n}\right)$ and any coefficients $\left(n_{j}\right) \subset \mathbb{C},$ we have
		\begin{align}\label{d}
	\left\|\sum_{j} n_{j}\left|e^{ i t \mathcal{L} } u_{j}\right|^{2}\right\|_{ L^{\frac{n+1}{n}}\left(\mathbb{T} \times \mathbb{C}^{n}\right)} \leq C \left(\sum_{j}\left|n_{j}\right|^{\frac{2(n+1)}{2n+1}}\right)^{\frac{(2n+1)}{2 (n+1)}},
		\end{align}
	where $C$ is a  constant depends on $n$.
\end{thm}
%\begin{thm}\label{diag} 	Let 	$n \geq 1$ and let $ S=\{(\mu, \nu, \lambda)\in \mathbb{N}_0^n\times \mathbb{N}_0\times \mathbb{Z}: \lambda=2|\nu|+n\}$ be the  discrete  surface defined. Then	$$\left\|W_{1} T_{S} W_{2}\right\|_{\mathcal{G}^{2(n+1)}\left(L^{2}\left(\mathbb{T}\times \mathbb{C}^n\right)\right)}  \leqslant C\left\|W_{1}\right\|_{L_{t, x}^{2(n+1)}\left(\mathbb{T}\times \mathbb{C}^n\right)} \left\|W_{2}\right\|_{L_{t, x}^{2(n+1)}\left(\mathbb{T}\times \mathbb{C}^n\right)}$$ for all $W_1, W_2$ with a constant $C > 0$ independent of $W_1, W_2.$ \end{thm}
\begin{pf} In order to prove (\ref{d}), by  Lemma \ref{dual},  it is enough to show \begin{eqnarray}\label{last}\left\|W_{1} T_{S} W_{2}\right\|_{\mathcal{G}^{2(n+1)}\left(L^{2}\left(\mathbb{T}\times \mathbb{C}^n\right)\right)}  \leq C\left\|W_{1}\right\|_{L_{t, w}^{2(n+1)}\left(\mathbb{T}\times \mathbb{C}^n\right)} \left\|W_{2}\right\|_{L_{t, w}^{2(n+1)}\left(\mathbb{T}\times \mathbb{C}^n\right)}\end{eqnarray} for all $W_1, W_2\in L_{t, w}^{2(n+1)}\left(\mathbb{T}\times \mathbb{C}^n\right),$  where  $T_S:=\mathcal{E}_S(\mathcal{E}_S)^*$.
For $-1<{\rm  Re}~z\leq0$, define the operator $T_z$ (on Schwartz class functions on $\mathbb{T}\times \mathbb{C}^n$) by\begin{align}\label{*}
T_z g(t,  w)  &=  \sum_{\mu, \nu, \lambda} \hat{g}(\mu, \nu, \lambda) G_z( \mu, \nu, \lambda) \Phi_{\mu\nu}(w) e^{-i\lambda t },\end{align} where $G_z$ is defined in (\ref{s11}). When ${\rm  Re}~z=0$, we have
	\begin{align}\label{two}
	\|T_{is} \|_{L^2(\mathbb{T}\times \mathbb{C}^n)\to L^2(\mathbb{T}\times \mathbb{C}^n)} =\left\|G_{i s}\right\|_{L^{\infty}\left(\mathbb{T}\times \mathbb{C}^n \right)}\leq \left|  \frac{1}{\Gamma(1+i s)} \right|\leq C e^{\pi|s| / 2}.	\end{align}
Further, using (\ref{*}), we have \begin{align*}
T_z g(t,  w)
&=  \sum_{\mu, \nu, \lambda} \int_{\mathbb{T}}(\mathcal{F}_3^{-1}\hat{g})(\mu, \nu)(s) G_z( \mu, \nu, \lambda) \Phi_{\mu\nu}(w) e^{-i\lambda (t-s) }~ds,
\end{align*}
where $\mathcal{F}_3^{-1} \hat{g}$ denotes the inverse Fourier transform of $\hat{g}$ with respect to third variable. Then using (\ref{s11}), Proposition \ref{515}  and  the distributional calculus of $(\lambda-(2|\nu|+n))_+^z$, we get
\begin{align}\label{K_z}\nonumber
T_z g(w,  t)
&=  \sum_{\mu, \nu, \lambda} \int_{\mathbb{T}}\langle {g}(\cdot, \cdot, s), \Phi_{\mu\nu}\rangle \Phi_{\mu\nu}(w)G_z( \mu, \nu, \lambda) e^{-i\lambda (t-s) }~ds\\\nonumber
&=  \sum_{\mu, \nu} \int_{\mathbb{T}}\langle {g}(\cdot, \cdot, s), \Phi_{\mu\nu}\rangle \Phi_{\mu\nu}(w)\sum_{\lambda}G_z( \mu, \nu, \lambda) e^{-i\lambda (t-s) }~ds\\\nonumber
&=  \frac{1}{\Gamma(z+1)}\sum_{\mu, \nu} \int_{\mathbb{T}}\langle {g}(\cdot, \cdot, s), \Phi_{\mu\nu}\rangle \Phi_{\mu\nu}(w) e^{-i(t-s)(2|\nu|+n)}\sum_{\kappa=0}^\infty \kappa_+^z  e^{-i\kappa (t-s) }~ds\\\nonumber
&=  -ie^{-iz\frac{\pi}{2}} \sum_{\mu, \nu} \int_{[-\pi,\pi]}\langle {g}(s, \cdot, \cdot), \Phi_{\mu\nu}\rangle \Phi_{\mu\nu}(w)  (t-s)^{-z-1}e^{-i(t-s)(2|\nu|+n)} ds\\
&\qquad +\frac{b(t)}{\Gamma(z+1)} \sum_{\mu, \nu} \int_{[-\pi,\pi]}\langle {g}(s, \cdot, \cdot), \Phi_{\mu\nu}\rangle \Phi_{\mu\nu}(w) e^{-i(t-s)(2|\nu|+n)}\;ds.
\end{align}
We ignore the second term in (\ref{K_z}) as it vanishes as $z$ tends to $-1$ in our calculation.
Now
\begin{align*}
&\sum_{\mu, \nu}  \langle {g}(\cdot, \cdot, s), \Phi_{\mu\nu}\rangle \Phi_{\mu\nu}(w) e^{-i(t-s)(2|\nu|+n)}\\
&=(2\pi)^{-\frac{n}{2}}\sum_{ \nu}  {g}(\cdot, \cdot, s)\times  \Phi_{\nu\nu}(w) e^{-i(t-s)(2|\nu|+n)}\\
&=(2\pi)^{-{n}}\sum_{ k=0}^\infty  e^{-i(t-s)(2k+n)} {g}(\cdot, \cdot, s)\times  \phi_{k}(w) \\
&=(2\pi)^{-{n}} {g}(\cdot, \cdot, s)\times  \sum_{ k=0}^\infty  e^{-i(t-s)(2k+n)}\phi_{k}(w) .
\end{align*}
Thus from (\ref{K_z}), we have
\begin{align*}
T_z g(t, w)
&=  -i(2\pi)^{-{n}} e^{-iz\frac{\pi}{2}}    \int_{\mathbb{T}}  {g}(\cdot, \cdot, s)\times  (t-s)^{-z-1} \sum_{ k=0}^\infty  e^{-i(t-s)(2k+n)}\phi_{k}(w) ds\\
&= -i(2\pi)^{-{n}} e^{-iz\frac{\pi}{2}}     \int_{\mathbb{T}}   \int_{\mathbb{C}^n}  {g}(u, s)H(u, w, t-s) e^{-\frac{i}{2} \operatorname{Im} (u \cdot \bar{w})} ds du\\
&=  -i(2\pi)^{-{n}} e^{-iz\frac{\pi}{2}}       \int_{\mathbb{C}^n}  ({g}(u, \cdot )*H(u, w, \cdot))(t) e^{-\frac{i}{2} \operatorname{Im} (u \cdot \bar{w})}   du,
\end{align*}
where $$H(u, w, t-s) =(t-s)^{-z-1} \sum_{ k=0}^\infty  e^{-i(t-s)(2k+n)}\phi_{k}(u-w).$$ When $z=-\lambda_0+is$, using (\ref{kk}), we get
\begin{align*}
|T_z g(t, w) |
&\leq  (2\pi)^{-{n}} e^{\frac{|s|\pi}{2}}      \int_{\mathbb{C}^n}  \|{g}(u, \cdot )\|_1\sup_{t\in [-\pi, \pi]}|H(u, w, t)|  du.
\end{align*}
Thus  $T_z$  is bounded from $ L^1(\mathbb{T}\times \mathbb{C}^n)$ to $ L^\infty (\mathbb{T}\times \mathbb{C}^n)$ if and only if $\sup_{t\in [-\pi, \pi]}|H(u, w, t)|$ is bounded for each $u, w\in\mathbb{C}^n.$ But from (\ref{kk}), we get
\begin{align}\label{kernel}
|H(u, w, t)|&\sim\frac{C}{ |t|^{{\rm Re}~(z+1+n)}} .
\end{align}
So for each  $w, u\in\mathbb{C}^n,$ $|H(u, w, t)|$ is bounded if and only if $\operatorname{Re}(z)=-(n+1).$   Thus \begin{eqnarray}\label{**}\|T_{-\lambda_0+is}\|_{L^1(\mathbb{T}\times \mathbb{C}^n)\to L^\infty (\mathbb{T}\times \mathbb{C}^n) }\leq C e^{-\frac{s\pi}{2}}.\end{eqnarray} By (\ref{two}) and (\ref{**}), the family of operators $(T_z)$  satisfy (\ref{expo}) and (\ref{expo1}). The conclusion of the theorem follows  by choosing $\lambda_0=n+1$  in Proposition \ref{sh4} and the identity $T_S=T_{-1}.$
\end{pf}
\subsection{Proof of  Theorem \ref{STHH}}
Using the fact that the operator $e^{-i t \mathcal{L}}$ is unitary, triangle inequality gives (\ref{SH1}) for the pair $(p,q)=(\infty,1)$. Equivalently, the operator $$W\in L_t^\infty L_z^2(\mathbb{T}\times\mathbb{C}^n)\mapsto We^{-itH}(e^{-itH})^*\overline{W}\in\mathcal{G}^\infty $$ is bounded by Lemma \ref{dual}.
Similarly, by (\ref{last}), the operator $$W\in L_t^{2(n+1)} L_z^{2(n+1)}(\mathbb{T}\times\mathbb{C}^n)\mapsto We^{-itH}(e^{-itH})^*\overline{W}\in\mathcal{G}^{2(n+1)} $$ is bounded.
Applying  the complex interpolation method \cite{ber} (Chapter  4), the operator $$W\in L_t^{\frac{2q}{2-q}} L_z^{\frac{2p}{2-p}}(\mathbb{T}\times\mathbb{C}^n)\mapsto We^{-itH}(e^{-itH})^*\overline{W}\in\mathcal{G}^\alpha $$ is bounded for $2\leq \frac{2q}{2-q}\leq 2(n+1)$ and $2(n+1)\leq \frac{2p}{2-p}\leq\infty$.  Again, applying Lemma \ref{dual}, the inequality (\ref{SH1}) holds for the range $1\leq q\leq 1+\frac{1}{n}$.

%When $q=1$ and $p=\infty$, using triangle inequality and the fact that $e^{it\mathcal{L}}$ is unitary operator on  $L^2(\mathbb{C}^n)$ for any fixted $t\in [-\pi, \pi]$, we have the bound\begin{align}\label{1}\sup _{t \in [-\pi, \pi ]}\left(\int_{\mathbb{C}^{n}}\left|\sum_{j} n_{j}\left | \left(e^{-i t \mathcal{L}} u_{j}\right) (x) \right |^{2} \right | d x\right) \leqslant \sum_{j}\left|n_{j}\right|.\end{align}  Further by Theorem  \ref{diag}  we have that (\ref{SH1}) holds for   $ p=q=1+\frac{1}{n} $. Applying  the complex interpolation method \cite{ber}(chapter  4) we  have Theorem \ref{STHH}.

\section*{Acknowledgments}
The first author thanks the Ministry of Human Resource Development, India for the  research fellowship and Indian Institute of Technology Guwahati for the support provided during the period of this work.

\end{document}